# FLUID-STRUCTURE INTERACTION PROBLEMS: AN APPLICATION TO ANCHORED AND UNANCHORED STEEL STORAGE TANKS SUBJECTED TO SEISMIC LOADINGS

Hoang Nam PHAN[1], Fabrizio PAOLACCI[2]

## ABSTRACT

The estimation of the degree of risk of steel storage tanks from an industrial or nuclear power plant during an earthquake event is a very difficult task since the dynamic behaviour of the tank-liquid system is highly complex. The seismic response of steel storage tanks is quite different from building structures not only due to hydrodynamic effects acting on the shell and bottom plate but also because of many sources of nonlinear behaviour mechanisms. For the earthquake-resistant design of tanks, it is important to use a rational and reliable nonlinear dynamic analysis procedure, which is capable of capturing the seismic behaviour of the tanks under artificial or real earthquakes. The present paper deals with the nonlinear finite element modelling of steel storage tanks subjected to seismic loadings. The interaction effects of fluid and structure are modelled using a surface-based acoustic-structural interaction in the ABAQUS software. The successive contact and separation between the bottom plate and its rigid foundation caused by sliding and uplift are taken into account by a contact modelling approach. Reduced scales of both slender and broad cylindrical steel storage tanks from a shaking table campaign within the framework of the INDUSE-2-SAFETY project are selected for the study. The seismic responses of the tanks including the hydrodynamic pressure distribution acting on the shell plate, the elevation of the liquid free surface, and the uplift response of the bottom plate when the tank is unanchored are presented and compared with the experimental results. The responses obtained from the refined models are also compared with those obtained using a simplified approach, which is based on the lumped mass model of the liquid motion.

## 1. INTRODUCTION

Industrial plants are complex systems of numerous integrated components and processes that can make them particularly vulnerable to natural hazard events. Due to the interaction between the natural events and the industrial risk, several effects take place in industrial plants, in particular, storage sites, causing damage to pipelines, process equipment, storage tanks, and consequently the release of hazardous materials. Earthquake damage in recent decades (e.g., 1995 Hyogoken-Nanbu Japan, 1999 İzmit Turkey, 2003 Tokachi-oki Japan, 2008 Wenchuan China, 2011 Great East Japan, 2012 Emilia Italy, etc.) has revealed that storage tanks are one of the most vulnerable components in industrial plants. Damage to tanks can cause significant disruption to the facility operation. Seriously, the extensive seismic-induced uncontrolled fires, when flammable materials or hazardous chemicals leak, naturally increase the overall damage to nearby areas.

The seismic response of storage tanks is quite different from buildings not only due to hydrodynamic effects acting on the shell plate but also because of many sources of nonlinear behaviour mechanisms. These effects may result in different failure modes of tank components during earthquakes, e.g., elastoplastic (elephant's foot) and elastic (diamond shape) buckling of the shell plate, rupture of the shell-to-bottom connection, damage of the fixed or floating roof, failure of the support structure, etc.

---

[1]PhD student, Roma Tre University, Rome, Italy, hoangnam.phan@uniroma3.it
[2]Assistant professor, Roma Tre University, Rome, Italy, fabrizio.paolacci@uniroma3.it

There has been an increasing trend toward the modelling of steel storage tanks subjected to seismic loadings. The analysis procedure has been commonly based on three possible models: (i) the simplified spring-mass model in which the impulsive and convective components are modelled as SDOF systems, (ii) the added mass model in which the impulsive and convective forces are converted to equivalent masses along the height of the shell and bottom plate, and (iii) the full nonlinear finite element (FE) model in which the real interaction between fluid and structure is considered. The basic idea of the simplified spring-mass model is based mainly on the spring-mounted masses analogy proposed by Housner (1963). This analogy is derived from the solution of the hydrodynamic equations that describe the behaviour of liquid inside a rigid container, fixed to the foundation. The solution of the impulsive hydrodynamic pressure, in some cases, is used as distributed masses on the shell plate, which is modelled using shell elements (Virella et al. 2006). The added mass method could be used in anchored tanks; however, it failed to model unanchored tanks because the uplift mechanism significantly affects the hydrodynamic pressure distribution of the tanks. This mechanism, which is developed in response to large overturning moment, is the dominant response under seismic loads. Associated with the base uplift, significant plastic rotations of the shell-to-bottom plate joint and intensive stresses in the shell plate are developed. Malhotra and Veletsos (1994) proposed a spring-mass model for the unanchored tanks, where the uplift behaviour of the bottom plate is modelled by idealising the bottom plate as a series of uniformly loaded semi-infinite, prismatic beams that rest on a rigid foundation. The approximate nature of the methodology suggested that extensive researches on the uplift behaviour might need a more refined approach, which is based on the nonlinear FE model of the tank-liquid system.

The numerical techniques, especially the finite element method has become a highly useful tool; it has possibly been used for the numerical analysis of not only the tank itself but also the contained liquid with more reliable analysis results. However, because of the complex nonlinear behaviour of the coupled tank-liquid-foundation system, the modelling of this type of structure is a very challenging topic, and many studies are still being carried out in this field. The behaviour of the liquid is commonly represented by the Lagrangian approach, where the liquid is assumed to be linearly elastic, inviscid, and irrotational (Bayraktar et al. 2010, Phan et al. 2017). This approach can be used to model the liquid motion in a rigid or flexible container; however, mesh-based Lagrangian methods have limitations with large deformations. If the liquid undergoes large deformations, the mesh has to be restructured to accommodate the new configuration. In this case, an Arbitrary Lagrangian-Eulerian (ALE) adaptive mesh may be used in the liquid domain to permit large deformations, especially at the liquid free surface (Phan et al. 2017). A coupling approach, so-called Eulerian-Lagrangian coupling (CEL), is useful for tank sloshing simulations. The CEL allows for the interaction between the Lagrangian tank domain where the material is fixed to the mesh and the Eulerian fluid domain where the material can flow through the mesh. The use of Eulerian elements eliminates the problem of extreme element deformation associated with Lagrangian fluid meshes (Tippmann et al. 2009, Mittal et al. 2014). However, this approach faces the challenge in modelling the contact between the fluid domain and the steel tank. When a thin shell plate is considered, the leakage of the material may occur at the interface during the analysis. Moreover, this approach is quite time-consuming because a more refined mesh needs to be used in the interface to avoid any penetration of the fluid. Therefore, a simplification of the FE model is commonly used through considering the liquid as inviscid, irrotational, and with no mean flow. This leads to using the acoustic wave equation which considers the propagation of the vibrating waves inside the fluid. Thus, this is the simplest formulation to take into account the fluid-structure interaction and thus less computational resources are required. This method is called as "structural-acoustic coupling" and adopted in the ABAQUS software (SIMULIA 2014).

The present study deals with the nonlinear FE modelling of slender and broad steel storage tanks subjected to seismic loadings. The modelling approach is capable of handling all complexities related to nonlinear behaviour mechanisms of the tank-liquid system. The numerical modelling is developed using the ABAQUS software with an implicit dynamic analysis. In particular, the shell and bottom plate are modelled using Lagrangian shell elements, while an acoustic FE mesh is used in the liquid domain. The tank-liquid interaction is considered using a surface-based tie constraint between the inner surface of the tank and the liquid domain. The successive contact and separation between the tank bottom plate and its rigid foundation caused by sliding and uplift of the bottom plate are also taken into account for the broad tank by a contact modelling approach. A comparative study of the seismic response of both slender and broad tanks is presented. The results in terms of the hydrodynamic pressure distribution, the maximum sloshing of the liquid free surface, and the uplift response of the bottom plate obtained from the FE models are evaluated and compared with those obtained from experiments and simplified models.

## 2. NUMERICAL MODELLING

### 2.1 Description of case studies

To facilitate the development and verification of the FE models, the sample tanks chosen for this study are based on a shaking table campaign within the framework of the INDUSE-2-SAFETY project (INDUSE-2-SAFETY 2016). The tests on the slender and broad tanks were conducted at the CEA Seismic Mechanic Studies Laboratory (EMSI) in Saclay (France).

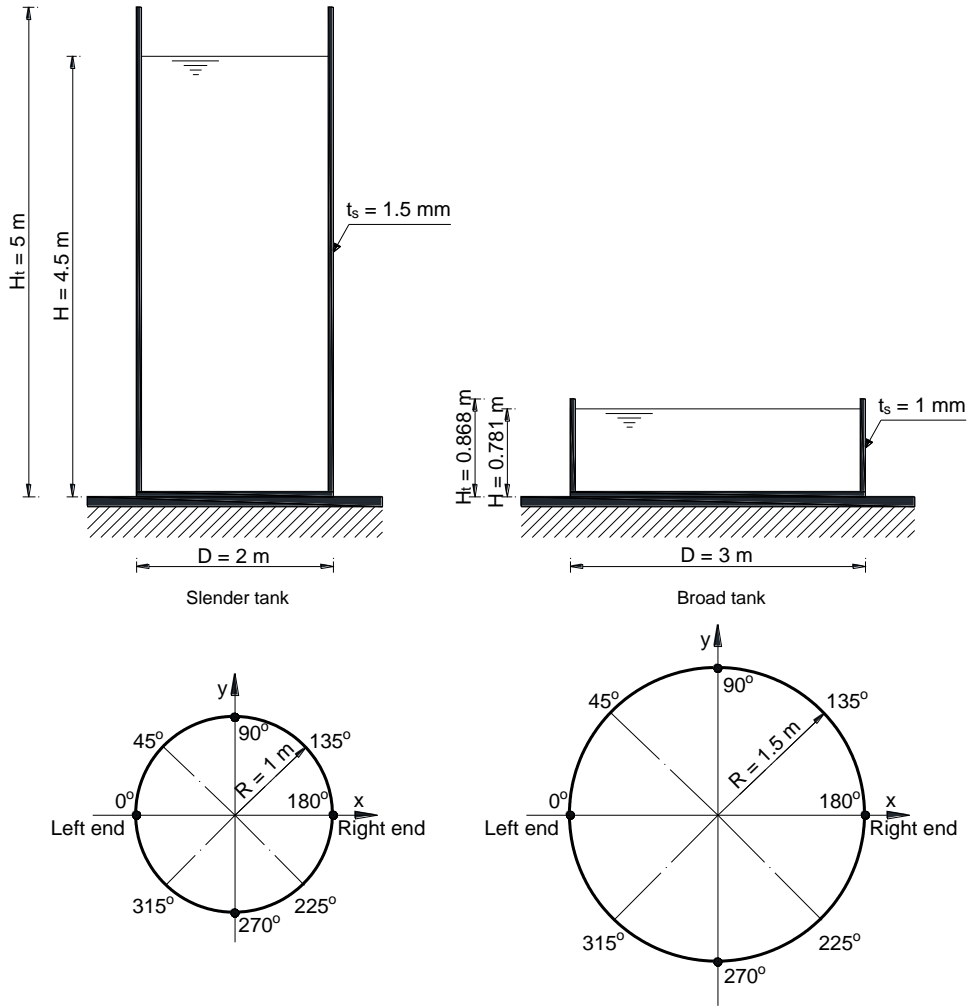

Figure 1. Schematic of the examined slender and broad tanks

The slender tank has a 2-m-diameter and a 5-m-height, which is scaled from a real tank with a scale ratio of 1/4, as shown in Figure 1. It is formed by the cylindrical S355JR steel sheet. The shell plate thickness is 1.5 mm. The base is welded to a massive basis, and a circular stiffener is welded to the top of the cylindrical shell. The mass of the empty tank is 2300 kg. Once filled with water at 90% of its height (i.e., 4.5 m), the specimen mass is 16.4 tons. The broad tank, having a diameter of 3 m and a total height of 0.868 m, is scaled from an existing broad tank at a refinery in Italy with a scale ratio of 1/18. The schematic of the reduced-scale tank, which includes both vertical and horizontal planes, is shown in Figure 1. The tank shell is formed by a cylindrical SS304 stainless steel sheet that has a 1-mm thick. The shell plate is welded to a bottom plate of the same material and thickness. The top of the shell is reinforced with a circular stiffener. The specimen is positioned in the centre of the table plate on an intermediate EPDM sheet. This rubber membrane helps increase friction coefficient between the bottom plate and the table to decrease the sliding effect and protect the mechanical bearings and electrical circuits in case of water overtopping. The estimated mass of the empty tank is 123 kg. The tank is filled with water at 90% of its height (i.e., 0.781 m), resulting in a total mass of 5.6 tons.

*2.2 Nonlinear FE modelling*

The nonlinear dynamic analysis of the tank-liquid system is carried out using the FE analysis package ABAQUS (SIMULIA 2014). The coupled acoustic-structural analysis is used in this study. This approach is simple and effective to treat numerically, as it assumes no material flow and thus no mesh distortion. The FE meshes of steel tank consist of four-node, doubly curved quadrilateral shell elements (S4R). Each node of shell element has three translational and three rotational degrees of freedom. The liquid is modelled using eight-node brick acoustic elements (AC3D8). The acoustic FE model is based on the linear wave theory and considers the dilatational motion of the liquid. To derive the equations for acoustic wave propagation, a number of assumptions have to be made to simplify the equations of fluid dynamics. The acoustic element has only one pressure unknown as the degree of freedom at each node. Hence, no actual flow occurs in an acoustic simulation. The tank-liquid interaction is considered using a surface-based tie constraint between the tank inner and liquid surface. This constraint is formulated based on a master-slave contact method, in which normal force is transmitted using tied normal contact between both surfaces through the simulation. When the system is subjected to large deformations, the ALE framework can be used to prescribe the movement of the acoustic mesh, including interior nodes, to follow and adapt to the movement of the structure. The sloshing waves are considered in the liquid model. Assuming the small-amplitude gravity waves on the liquid free surface, the pressure boundary condition specified at the free liquid surface can be presented in the form of Equation 1, where $p$ is the hydrodynamic pressure at the liquid free surface.

$$\frac{\partial^2 p}{\partial t^2} + g\frac{\partial p}{\partial z} = 0 \qquad (1)$$

The bottom plate of the slender tank is assumed to be rigid and fixed to the foundation, while the broad one is unanchored and rested on a rigid slab that is modelled using solid elements. The successive contact and separation between the tank bottom plate and its rigid foundation are taken into account by a surface-based contact modelling algorithm. The boundary conditions of the model are shown in Figure 2.

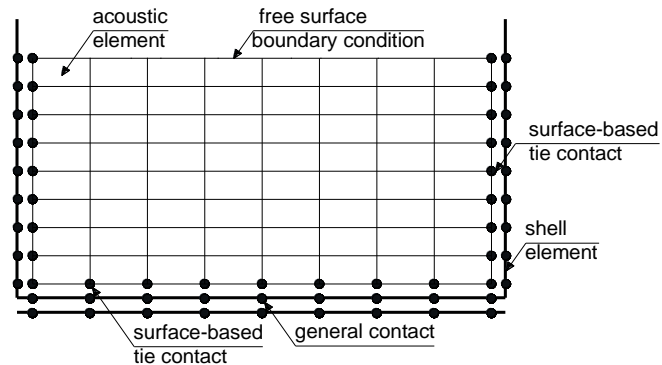

Figure 2. Boundary conditions of the liquid-tank model

Both geometric and material nonlinearities are considered in the analysis. The plasticity of the steel tank is modelled based on the stress-strain curve of the material. The curve obtained from mechanical testing is converted into the true stress and the plastic strain. The water density is considered to be 998.21 kg/m$^3$, and its bulk modulus is 2150 MPa. The Rayleigh mass proportional damping is employed for the tank model assuming a damping ratio of 2.0%, for the fundamental vibration mode of the tank-liquid system. Due to the structural symmetry and to reduce the computational cost, only half of the tank-liquid system is modelled, and symmetry plane boundary conditions are employed. The mesh convergence analysis results in an optimal mesh size of 0.04 m and 0.08 m in the longitudinal direction and the circumferential direction, respectively, to achieve acceptable accuracy. An example of the FE mesh of the broad tank model is illustrated in Figure 3.

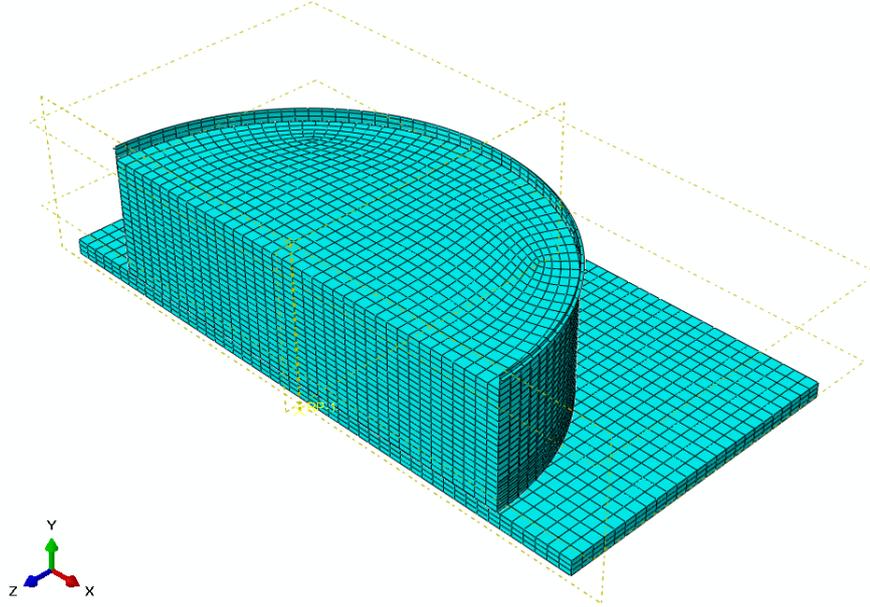

Figure 3. Numerical model of the tank-liquid system

The acoustic wave equations do not include any terms for body forces, which means that forces such as gravity are not included. Hence, at the first step, the tank is subjected to the gravity load and the hydrostatic pressure acting on the shell and bottom plate. The hydrodynamic pressures acting on the shell and bottom plate are measured during the dynamic analysis through the tie contact between the shell and acoustic elements.

*2.3 Modal analysis*

The modal analyses of the two tanks are first performed. The natural frequency calculation is based on the Lanczos eigensolver method. The acoustic-structural coupling is projected onto the subspace of eigenvectors using SIM-based linear dynamic procedures. The natural periods of the two tanks computed by the modal analyses are in close agreement with those obtained from Eurocode 8 (EN 1998-4 2006), as shown in Table 1.

Table 1. Modal analysis results

| Vibration mode | Slender tank | | Broad tank | |
|---|---|---|---|---|
| | **FE model** | **EN 1998-4** | **FE model** | **EN 1998-4** |
| First impulsive, $T_i$ (s) | 0.061 | 0.069 | 0.013 | 0.016 |
| First convective, $T_{c1}$ (s) | 1.478 | 1.479 | 2.100 | 2.100 |
| Second convective, $T_{c2}$ (s) | 0.867 | 0.869 | 1.068 | 1.068 |
| Third convective, $T_{c3}$ (s) | 0.679 | 0.687 | 0.895 | 0.841 |

*2.4 Description of input signals*

The horizontal components of the L'Aquila Italy 06/04/2009, the Chi-Chi Taiwan 21/09/1999, and Northridge USA 17/01/1994 earthquakes have been used. The seismic signals used for the tests were chosen among the ones maximising the mechanical effects on the tanks. In particular, the Chi-Chi signal is for the sloshing behaviour of the broad tank, the Northridge signal is for the uplift behaviour of the broad tank, and the L'Aquila signal is for the buckling behaviour of the slender tank.

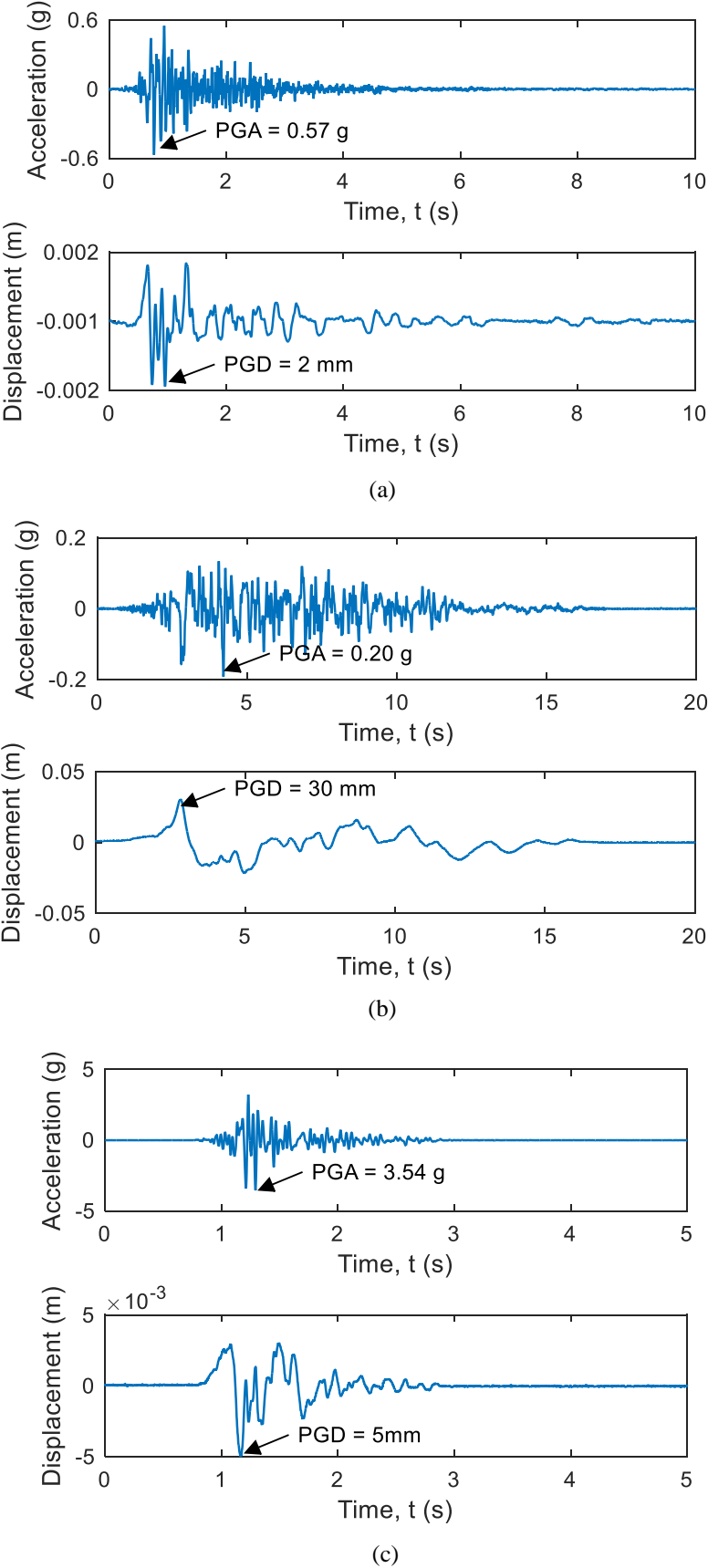

Figure 4. Input signals for the tests: (a) L'Aquila signal, (b) Chi-Chi signal, and (c) Northridge signal

Different time scales based on the Froude similarity are defined for each test, depending on the specific response of the tanks. The scaled signals for the tests are shown in Figure 4 in terms of acceleration and displacement time histories. The peak ground acceleration (PGA) and peak ground displacement (PGD) of each signal are also shown in the plots.

## 3. COMPARATIVE RESULTS

### 3.1 Slender tank

The seismic response analysis of the slender tank, fixed to the foundation, is first performed using the L'Aquila signal with a PGA level about 0.57 g. The analysis results are obtained and compared with the experimental data. A simplified approach, which was presented in Phan et al. (2016), is also used to calculate the seismic response of the tank. The approach was based mainly on the formulas in Eurocode 8 (EN 1998-4 2006). The results in terms of the hydrodynamic pressure distribution on the shell obtained from three models are shown in Figure 5.

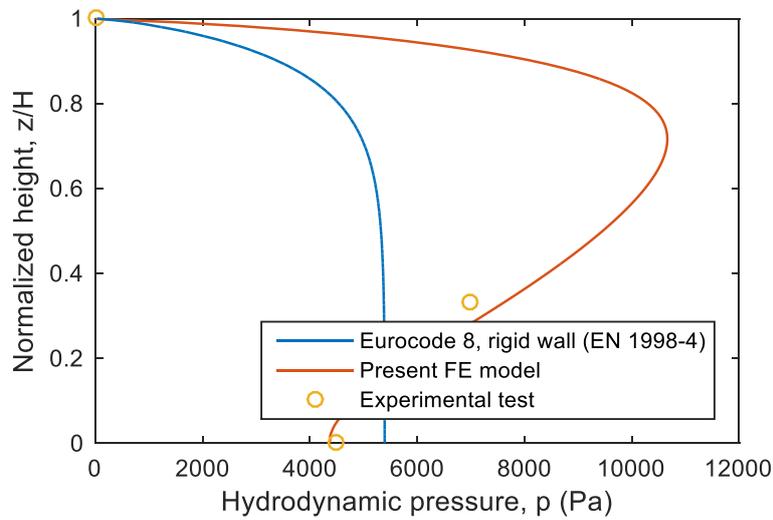

Figure 5. A comparison of peak hydrodynamic pressure distribution on the shell plate at $t = 0.947$ s

It can be observed that the pressure distribution obtained from the present FE model for the slender tank varies irregularly, exhibiting its maximum value away from the bottom. This distribution shape is in contrast to those of the spring-mass model, where the pressure increases monotonically from top to bottom of the tank. The hydrodynamic pressure values measured from the experimental test at some positions on the shell plate are also shown in the figure. These pressures are in good agreement with the ones of the present model. The pressure obtained using the formulas in Eurocode 8 (EN 1998-4 2006) is slight higher than that of the refined model measured at the bottom; however, from $z/H = 0.2$ to the top, the pressure distribution obtained using the code is significantly smaller than of the present model. It can be concluded that an accurate solution of the hydrodynamic pressure that is obtained for a flexible shell plate is needed when slender tanks are considered.

### 3.2 Broad tank

The hydrodynamic pressure on the shell plate is obtained for the broad tank from the nonlinear time history dynamic analysis, which uses the Chi-Chi signal. Figure 6 shows the comparison between three results of the hydrodynamic pressure distribution along the tank height that are obtained from the spring-mass model, present FE model, and experimental test. It is noticed that during a strong seismic action, the uplift of the bottom plate may occur. A rotation spring representing the uplift resistance of the bottom plate is added to the spring-mass model, as presented in Phan et al. (2016). The figure shows a good agreement between the results obtained from the refined model and the experimental test. Differently from the slender tank, the distribution shapes of the pressure for three cases are quite similar, where the pressures increase monotonically from top to bottom of the tank. However, the amplitude of the pressure obtained from the spring-mass model has a slight difference as compared with the ones from the refined model and the test. In particular, a considerable increase is recorded at locations on the shell plate near the top. The time history responses of the hydrodynamic pressure at the left end of the bottom plate are shown in Figure 7.

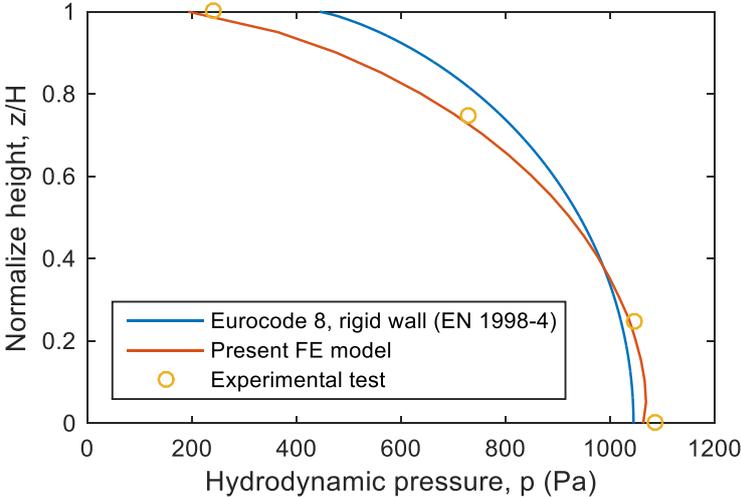

Figure 6. A comparison of peak hydrodynamic pressure distribution on the shell plate at $t = 4.04$ s

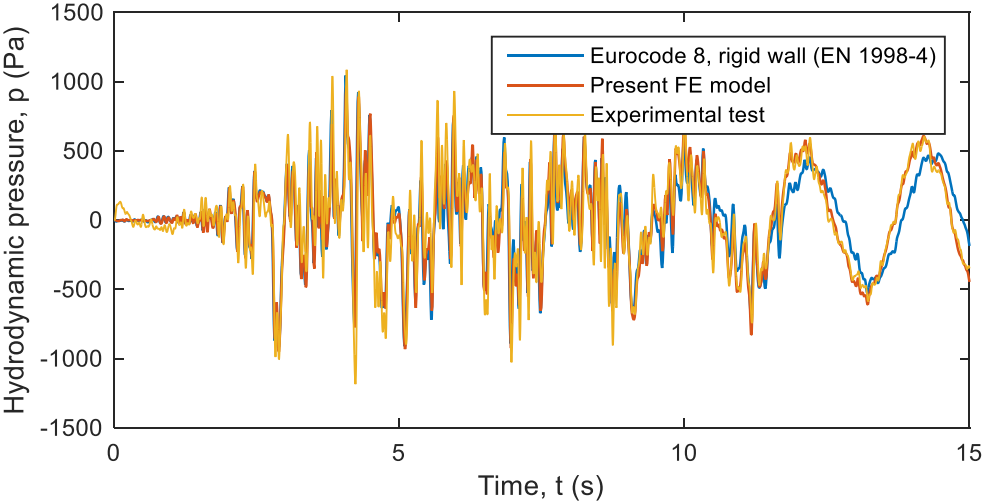

Figure 7. A comparison of hydrodynamic pressure time history at the left end of the bottom plate

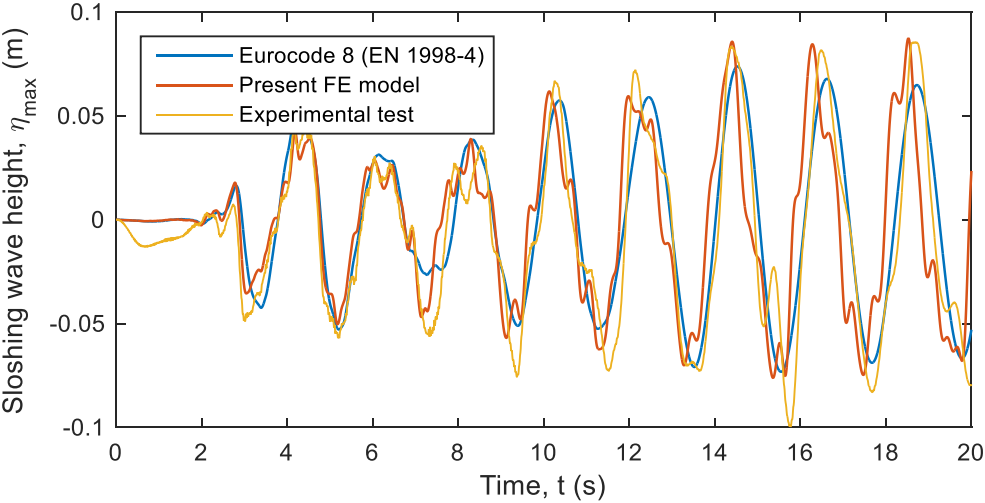

Figure 8. A comparison of maximum sloshing wave height time history at the left end of the liquid free surface

The time history responses of the sloshing wave height at the left end of the free surface are presented in Figure 8. The results are obtained using the Chi-Chi signal. It can be seen from the figure that the numerical results are highly consistent with those obtained from the experimental test in terms of both frequency and amplitude. The maximum sloshing wave heights measured from the present model and the experimental data are summarised in Table 2. The experimental data observed is slightly higher than the numerical ones. It is also noticed from the experimental data that the maximum sloshing wave height measured is equal to the freeboard height of the tank. An overtopping of the contained liquid was observed in the test.

Figure 9 presents the time history responses of the base uplift of the tank model measured at the left end of the bottom plate, and their peak values are shown in Table 2. The responses are obtained with the Northridge input signal. The uplift displacements of the experimental study include negative values because the tank is settled on an EPDM rubber. The numerical and experimental models lead to a relatively accurate description of the base uplift displacement at the right end for the input earthquake motion; however, there are some parts of the history response at the left end showing slight differences among three obtained results. In general, it can be concluded that the spring-mass model for the broad tank is capable of capturing the peak responses compared with the detailed FE analysis and the test.

Table 2. Peak response results of the broad tank

| Demand | Spring-mass model | FE model | Experimental test |
|---|---|---|---|
| Sloshing wave height, $\eta_{max}$ (m) | 0.074 | 0.081 | 0.086 |
| Uplift displacement, $w$ (m) | 0.016 | 0.020 | 0.021 |

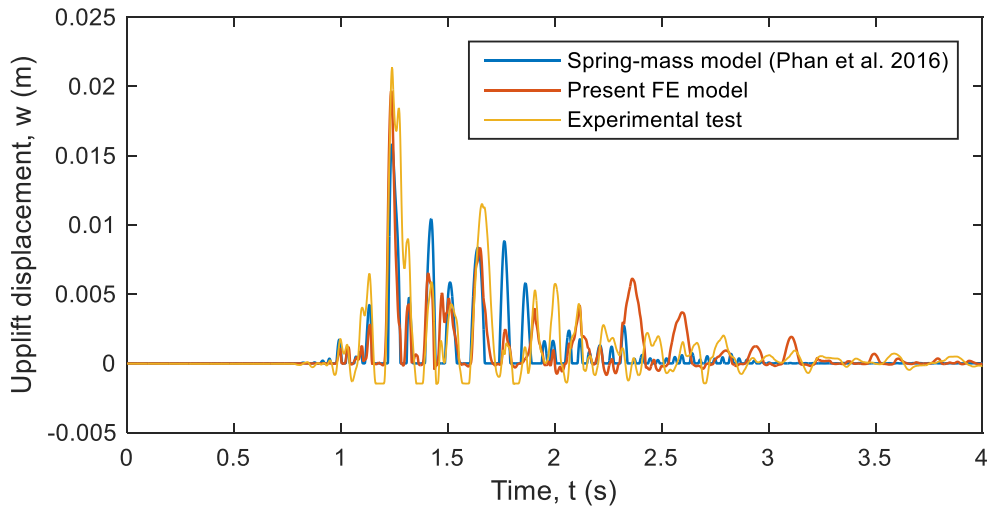

Figure 9. A comparison of uplift time history at the left end of the bottom plate

## 4. CONCLUSIONS

In this study, a FE modelling approach has been presented for the seismic response analysis of both anchored and unanchored steel storage tanks. The approach is based on the structural-acoustic coupling and is implemented in the ABAQUS software. The complex interaction mechanisms between the shell plate and the liquid, as well as the bottom plate and the rigid foundation, have been taken into account with contact algorithms. As case studies, reduced-scale slender and broad tanks, which were used in shaking table tests, have been investigated in the study. The time history dynamic analyses have been performed on the tank models using three input signals that were selected for the shaking table tests. The records were scaled based on a suitable similarity approach to maximise the seismic responses of the tanks. The analysis results in terms of the hydrodynamic pressure distribution, the free surface sloshing, and the base uplift have been obtained and compared well with those observed from the experimental tests. The results of simplified models for the two case studies have also been presented to confirm the possibility of using the formulas in current design codes. It can be concluded from the compatible results between the numerical and experimental studies that the present FE model is capable of producing reliable seismic responses of steel storage tanks.

## 5. ACKNOWLEDGMENTS

This work has been partially funded by the European project INDUSE-2-SAFETY (Grant No. RFS-PR13056). The authors acknowledge CEA Seismic Mechanic Studies Laboratory who provided shaking table test data of the tanks.

## 6. REFERENCES

Bayraktar A, Sevim B, Altunışık AC, Türker T (2010). Effect of the model updating on the earthquake behavior of steel storage tanks. *Journal of Constructional Steel Research*, 66(3): 462-469.
Housner GW (1963). The dynamic behaviour of water tanks. *Bulletin of the Seismological Society of America*, 53: 381-387.
INDUSE-2-SAFETY (2016). Component fragility evaluation, seismic safety assessment and design of petrochemical plants under design-basis and beyond-design basis accident conditions. Annual report from 01/01/2016 to 31/12/2016, University of Trento, Italy.
Malhotra PK, Veletsos AS (1994). Uplifting response of unanchored liquid-storage tanks. *Journal of Structural Engineering*, 120(12): 3524-3546.
Mittal V, Chakraborty T, Matsagar V (2014). Dynamic analysis of liquid storage tank under blast using coupled Euler-Lagrange formulation. *Thin-Walled Structures*, 84: 91-111.
Ozdemir Z, Souli M, Fahjan YM (2010). Application of nonlinear fluid-structure interaction methods to seismic analysis of anchored and unanchored tanks. *Engineering Structures*, 32(2): 409-423.
Phan H, Paolacci F, Alessandri S (2016). Fragility analysis methods for steel storage tanks in seismic prone areas. *ASME Pressure Vessels and Piping Conference*, 8: V008T08A023.
Phan H, Paolacci F, Mongabure P (2017). Nonlinear finite element analysis of unanchored steel liquid storage tanks subjected to seismic loadings. *ASME Pressure Vessels and Piping Conference*, 8: V008T08A040.
SIMULIA (2014). Abaqus 6.14 Documentation. Dassault Systèmes Simulia Corp., Providence, RI, USA.
Tippmann J, Prasad S, Shah P (2009). 2-D tank sloshing using the coupled Eulerian-Lagrangian (CEL) capability of ABAQUS®/Explicit. *2009 SIMULIA Customer Conference*, Dassault Systemes Simulia Corp.
Virella JC, Godoy LA, Suárez LE (2006). Dynamic buckling of anchored steel tanks subjected to horizontal earthquake excitation. *Journal of Constructional Steel Research*, 62(6): 521-531.